\documentclass[11pt]{article}
\usepackage{mathrsfs}
\usepackage{amssymb}
\usepackage{amsfonts}
\usepackage{color,xcolor}
\usepackage{graphicx}
\usepackage{geometry}
\usepackage{manfnt}
\usepackage[pagewise]{lineno}
\usepackage[colorlinks,citecolor=blue,urlcolor=blue]{hyperref}
\newtheorem{theorem}{Theorem}[section]
\newtheorem{corollary}{Corollary}[section]
\newtheorem{lemma}{Lemma}[section]
\newtheorem{proposition}{Proposition}[section]

\newtheorem{defi}{Definition}[section]
\newtheorem{rem}{Remark}[section]

\def\[{{\Big[}}\def\]{{\Big]}}\def\({{\Big(}}\def\){{\Big)}}

\def\cD{{\mathcal D}}
\def\cF{{\mathcal F}}

\def\mE{{\mathbb E}}

\def\mN{{\mathbb N}}\def\mP{{\mathbb P}}
\def\mR{{\mathbb R}}

\def\={&\!\!=\!\!&}
\def\geq{\geqslant}\def\leq{\leqslant}

\begin{document}
\title{\bf A Kolmogorov type theorem for stochastic fields}

\author{Jinlong Wei$^a$, Guangying Lv$^{b}$\\
\\ {\small \it $^a$School of Statistics and Mathematics, Zhongnan University of}\\  {\small \it Economics and Law, Wuhan 430073, China} \\ {\small \tt  weijinlong.hust@gmail.com}\\ {\small \it $^b$College of Mathematics and Statistics, Nanjing University of Information} \\ {\small \it Science and Technology, Nanjing 210044, China}\\
{\small \tt gylvmaths@henu.edu.cn}}

\date{}
\maketitle

\noindent{\hrulefill}
\vskip1mm\noindent{\bf Abstract} We generalize the Kolmogorov continuity theorem and prove the continuity of a class of stochastic fields with the parameter. As an application, we derive the continuity of solutions for nonlocal stochastic parabolic equations driven by non-Gaussian L\'{e}vy noises.



\vskip1mm\noindent{\hrulefill}

\section{Introduction and main result}\setcounter{equation}{0}
In 1933, Kolmogorov \cite{Kol1933} proved the celebrate result as the following:
\begin{proposition}\label{pro1.1} If $\{X(t), \ t\in [0,1]\}$ is a stochastic process, such that
\begin{eqnarray}\label{1.1}
{\mathbb E}[|X(t)-X(s)|^\gamma]\leq C|t-s|^{1+\varepsilon},
\end{eqnarray}
where $\gamma$, $\varepsilon$ and $C$ are positive constants independent of $t$, then the trajectories of the process are continuous with probability 1.
\end{proposition}

The above result is one of the central aspects of stochastic analysis with many applications in the study of the asymptotic distributions of certain tests \cite{Kol1933}, the tightness of a sequence of stochastic processes \cite{Sko56,Sko57}, the existence of strong solutions and stochastic flows to stochastic differential equations \cite{Kuo}, and the regularization of local time for continuous semimartingale \cite{Yor}. Now, Proposition \ref{pro1.1} is well known for Kolmogorov's continuity theorem/criteron. Since then, Kolmogorov's result was strengthened in different forms \cite{Chen,Loe1,Loe2}. To be precise, we give a standard and simple type for the Kolmogorov theorem.
\begin{proposition}\label{pro1.2} \cite[Theorem 3.3.8,p31]{Kuo} and \cite[Theorem 2.1,p25]{RY}   Let $\{X(t), \ t\in [0,1]\}$ be a stochastic process defined on a probability space $(\Omega,\cF,\mP)$. Assume that there exist three strictly positive constants $\gamma,C$ and $\varepsilon$ such that for all $0\leq t,s\leq 1$, (\ref{1.1}) holds. Then $X$ has a continuous realization $\tilde{X}$, namely, there exists $\Omega_0$ such that $\mP(\Omega_0)=1$ and for each $\omega\in\Omega_0$, $X(t,\omega)=\tilde{X}(t,\omega)$ and $\tilde{X}(t,\omega)$ is a continuous function of $t$. Moreover, for every $\alpha\in [0,\varepsilon/\gamma)$,
\begin{eqnarray*}
{\mathbb E}\Big[\sup_{t\neq s}\Big(\frac{|\tilde{X}(t)-\tilde{X}(s)|}{|t-s|^\alpha}\Big)^\gamma\Big]<+\infty.
\end{eqnarray*}
\end{proposition}

In the present paper, we will generalize Proposition \ref{pro1.2} to the following form.

\begin{theorem} \label{the1.1} Let $(H,\ \|\ \|_H)$ be a Banach space and let $\{X_t(x), x\in [0,1]^d, t\in [0,1]\}$ be an $H$-valued stochastic field. Let $\varphi$ be a nonnegative and nondecreasing continuous function on
${\mathbb R}_+=(0,+\infty)$ such that $\lim_{r\rightarrow 0+}\varphi(r)=0$.  Suppose that there exist two strictly positive constants $\gamma$ and $C$ such that
\begin{eqnarray}\label{1.2}
{\mathbb E}[\sup_{0\leq t\leq 1}\|X_t(x)-X_t(y)\|_H^\gamma]\leq C|x-y|^d\varphi(|x-y|), \ \ x,y\in [0,1]^d,
\end{eqnarray}
and there exists another constant $0<\vartheta<1/\gamma$ such that when $\gamma\geq 1$,
\begin{eqnarray}\label{1.3}
\sum_{i=0}^\infty\varphi^\vartheta(2^{-i})<+\infty,
\end{eqnarray}
and when $\gamma<1$,
\begin{eqnarray}\label{1.4}
 \sum_{i=0}^\infty\varphi^{\gamma\vartheta}(2^{-i})<+\infty.
\end{eqnarray}
Further, we assume that there is a constant $\lambda\geq1$ such that for all sufficiently large natural number $n$,
\begin{eqnarray}\label{1.5}
\lambda^{-1}\leq \frac{\varphi(2^{-n})}{\varphi(2^{-n-1})}\leq \lambda.
\end{eqnarray}
Then there is a realization $\tilde{X}$ of $X$ such that $\tilde{X}$ is continuous in $x$ and for every $\alpha\in (0,1/\gamma-\vartheta]$
\begin{eqnarray}\label{1.6}
{\mathbb E}[\sup_{0\leq t\leq 1}[\tilde{X}_t(\cdot)]^\gamma_{\alpha,\varphi}]<+\infty,
\end{eqnarray}
where
\begin{eqnarray*}
[\tilde{X}_t(\cdot)]_{\alpha,\varphi}:=\sup_{x\neq y}\Big\{ \frac{\|\tilde{X}_t(x)-\tilde{X}_t(y)\|_H}{\varphi^\alpha(|x-y|)}, \ x,y\in [0,1]^d\Big\}.
\end{eqnarray*}
\end{theorem}

\begin{rem}\label{rem1.1} (i) Let $T$ and $c_1$ be two positive real numbers, and let $B_{c_1}$ be the ball centered at $0$ with radius $c_1$. By scaling transformations, the conclusions in Theorem \ref{the1.1} are true if one replaces $[0,1]$ and $[0,1]^d$ by $[0,T]$ and $B_{c_1}$, respectively.

(ii) Theorem \ref{the1.1} holds true as well if one uses a complete metric space $(S,\rho)$ instead of the Banach space $(H,\ \|\ \|_H)$.

(iii) When $\varphi(|x-y|)=|x-y|^\varepsilon$, then (\ref{1.5}) is true. Moreover, for every $\vartheta>0$,
\begin{eqnarray*}
\sum_{i=0}^\infty\varphi^\vartheta(2^{-i})=\sum_{i=0}^\infty   2^{-i\varepsilon\vartheta}<+\infty \quad \mbox{and} \quad \sum_{i=0}^\infty\varphi^{\gamma\vartheta}(2^{-i})=\sum_{i=0}^\infty   2^{-i\varepsilon\gamma\vartheta}<+\infty.
\end{eqnarray*}
Therefore, (\ref{1.3}) and (\ref{1.4}) hold. By Theorem \ref{the1.1}, there is a continuous realization $\tilde{X}$ of $X$ such that for every $\alpha\in [0,1/\gamma)$
\begin{eqnarray*}
{\mathbb E}\sup_{0\leq t\leq 1}
\sup_{x\neq y} \frac{\|\tilde{X}_t(x)-\tilde{X}_t(y)\|_H^\gamma}{|x-y|^{\varepsilon\alpha\gamma}}<+\infty.
\end{eqnarray*}
This is the result of Proposition \ref{pro1.2}, so we extend Kolmogorov's continuity theorem.
\end{rem}

Let $\gamma$ be given in Theorem \ref{the1.1} and let $\beta>1$ be a real number such that $\beta>\gamma$ if $\gamma \geq1$ and $\beta>1/\gamma$ if $\gamma<1$.  Set
\begin{eqnarray}\label{1.7}
\varphi(r)=\left\{
  \begin{array}{ll}
    (-\log(r))^{-\beta}, & \hbox{when} \quad r\leq \frac{1}{2}, \\
    \ \ \ f(r),& \hbox{when} \quad r>\frac12,
  \end{array}
\right.
\end{eqnarray}
where $f(r)$ is an arbitrary nondecreasing continuous function defined on $[1/2,+\infty)$ such that $f(\frac12)=(\log2)^{-\beta}$. If $\gamma\geq1$, we choose $1/\beta<\vartheta<1/\gamma$, then
\begin{eqnarray*}
\sum_{i=0}^\infty\varphi^\vartheta(2^{-i})=f^\vartheta(1)+
\sum_{i=1}^\infty(-\log(2^{-i}))^{-\beta\vartheta}=f^\vartheta(1)+
\frac{1}{(\log2)^{\beta\vartheta}}
\sum_{i=1}^\infty\frac{1}{i^{\beta\vartheta}}
<+\infty.
\end{eqnarray*}
If $\gamma< 1$, we fetch $1/(\beta\gamma)<\vartheta<1/\gamma$, then
\begin{eqnarray*}
\sum_{i=0}^\infty\varphi^{\gamma\vartheta}(2^{-i})=
f^{\gamma\vartheta}(1)+
\sum_{i=1}^\infty(-\log(2^{-i}))^{-\beta\gamma\vartheta}=f^{\gamma\vartheta}(1)+
\frac{1}{(\log2)^{\beta\gamma\vartheta}}
\sum_{i=1}^\infty\frac{1}{i^{\beta\gamma\vartheta}}
<+\infty.
\end{eqnarray*}
Further, for every $1\leq n\in\mN$,
\begin{eqnarray*}
\frac{\varphi(2^{-n})}{\varphi(2^{-n-1})}=\frac{(-\log(2^{-n}))^{-\beta}}{(-\log(2^{-n-1}))^{-\beta}}
=(1+\frac1n)^{\beta}.
\end{eqnarray*}
So (\ref{1.5}) holds true with $\lambda=2^\beta$.

Therefore, we have the following corollary.

\begin{corollary} \label{cor1.1} Let $(H,\ \|\ \|_H)$ be a Banach space and let $\{X_t(x), x\in [0,1]^d, t\in [0,1]\}$ be an $H$-valued stochastic field. Let $\gamma$ and $\beta$ be two real number such that $\beta>\gamma$ if $\gamma \geq1$ and $\beta>1/\gamma$ if $\gamma<1$. Suppose that (\ref{1.2}) holds with $\varphi$ given by (\ref{1.7}). Then there is a realization $\tilde{X}$ of $X$ such that $\tilde{X}$ is continuous in $x$ and for every $\alpha\in (0,1/\gamma-\vartheta]$, (\ref{1.6}) holds.
\end{corollary}

Let $\gamma$ and $\beta$ be given in Corollary \ref{cor1.1}. We choose $\varphi$ as the following form
 \begin{eqnarray}\label{1.8}
\varphi(r)=\left\{\begin{array}{ll}
 \frac{\log(-k_0\log(r))}{(-\log(r))^\beta}, & \hbox{when} \quad r\leq r_0, \\  \quad g(r),& \hbox{when} \quad r>r_0,
\end{array}\right.
\end{eqnarray}
where $k_0,r_0$ are two positive real numbers such that
$$
r_0=e^{-\frac{1}{k_0}e^{\frac{1}{\beta}}}.
$$
$g(r)$ is an arbitrary nondecreasing continuous function defined on $[r_0,+\infty)$ such that $g(r_0)=\log(-k_0\log(r_0))/(-\log(r_0))^\beta>0$. For $r\in (0,r_0)$, then
 \begin{eqnarray*}
\varphi^\prime(r)&=&
 \frac{1}{(-\log(r))^{2\beta}}\[(\log(-k_0\log(r)))^\prime(-\log(r))^{\beta}-
 \log(-k_0\log(r))((-\log(r))^{\beta})^\prime\] \\&=&\frac{1}{(-\log(r))^{2\beta}}\[-r^{-1}(-\log(r))^{\beta-1}+
 r^{-1}\beta\log(-k_0\log(r))((-\log(r))^{\beta-1})\] \\&=&\frac{1}{r(-\log(r))^{\beta+1}}\[\beta\log(-k_0\log(r))-1\]>0,
\end{eqnarray*}
so $\varphi$ is nonnegative and nondecreasing on $\mR_+$. Observing that $0<r_0<1$, there is a positive natural number $N_0$, such that $2^{-N_0}\leq r_0<2^{-N_0+1}$.

If $\gamma\geq1$, we choose $1/\beta<\vartheta<1/\gamma$, then
\begin{eqnarray*}
\sum_{i=0}^\infty\varphi^\vartheta(2^{-i})&=&\sum_{i=0}^{N_0-1}\varphi^\vartheta(2^{-i})+
\sum_{i=N_0}^\infty
\[ \frac{\log(-k_0\log(2^{-i}))}{(-\log(2^{-i}))^\beta}\]^\vartheta
\nonumber\\&\leq& \sum_{i=0}^{N_0-1}\varphi^\vartheta(2^{-i})+
\frac{1}{(\log2)^{\beta\vartheta}}\sum_{i=N_0}^\infty
\[ \Big(\frac{\log i)}{i^\beta}\Big)^\vartheta
+\Big(\frac{|\log(k_0\log2)|}{i^\beta}\Big)^\vartheta\]
\nonumber\\ &<&+\infty.
\end{eqnarray*}
If $\gamma< 1$, we fetch $1/(\beta\gamma)<\vartheta<1/\gamma$, then
\begin{eqnarray*}
\sum_{i=0}^\infty\varphi^{\gamma\vartheta}(2^{-i})&=&
\sum_{i=0}^{N_0-1}\varphi^{\gamma\vartheta}(2^{-i})+
\sum_{i=N_0}^\infty
\[ \frac{\log(-k_0\log(2^{-i}))}{(-\log(2^{-i}))^\beta}\]^{\gamma\vartheta}
\nonumber\\&\leq& \sum_{i=0}^{N_0-1}\varphi^{\gamma\vartheta}(2^{-i})+
\frac{1}{(\log2)^{\beta\gamma\vartheta}}
\sum_{i=N_0}^\infty
\[ \Big(\frac{\log i)}{i^\beta}\Big)^{\gamma\vartheta}
+\Big(\frac{|\log(k_0\log2)|}{i^\beta}\Big)^{\gamma\vartheta}\]
\nonumber\\ &<&+\infty.
\end{eqnarray*}
Moreover, for every $ \log_2(r_0^{-1})\leq n\in\mN$,
\begin{eqnarray*}
\frac{\varphi(2^{-n})}{\varphi(2^{-n-1})}=\frac{\log(-k_0\log(2^{-n}))}{(-\log(2^{-n}))^\beta}\times
\frac{(-\log(2^{-n-1}))^\beta}{\log(-k_0\log(2^{-n-1}))}
=(1+\frac1n)^{\beta}\frac{\log(k_0\log2)+\log(n)}{\log(k_0\log2)+\log(n+1)}.
\end{eqnarray*}
So (\ref{1.5}) holds true with $\lambda=\max\{2^\beta,
\frac{\log(k_0\log2(\log_2(r_0^{-1})+1))}{\log(k_0\log2\log_2(r_0^{-1}))}\}$.

Hence, we have the following result.
\begin{corollary} \label{cor1.2} Let $(H,\ \|\ \|_H)$, $\{X_t(x), x\in [0,1]^d, t\in [0,1]\}$, $\gamma$ and $\beta$ be stated in Corollary \ref{cor1.1}. Suppose that (\ref{1.2}) holds with $\varphi$ given by (\ref{1.8}). Then there is a realization $\tilde{X}$ of $X$ such that $\tilde{X}$ is continuous in $x$ and for every $\alpha\in (0,1/\gamma-\vartheta]$, (\ref{1.6}) holds.
\end{corollary}

We will give the proof details of Theorem \ref{the1.1} in Section 2, and as an application, we prove the continuity of solutions for nonlocal stochastic parabolic equations in Section 3.

\section{Proof of Theorem \ref{the1.1}}\setcounter{equation}{0}
Let $D_m$ be the set of points in $[0,1]^d$ whose components are equal to $2^{-m}i$ for some integral $i\in [0,2^m]$. Then $[0,1]^d=\cup_{m}D_m$. Let further $\Delta_m$ be the set of pairs $(x,y)$ in $D_m$ such that $|x-y|=2^{-m}$. There are no more than $d2^{(m+1)d}$ such pairs in $D_m$.

By setting  $K_i(t)=\sup_{(x,y)\in \Delta_i}\|X_t(x)-X_t(y)\|_H$, then there is a constant $C_1$ such that
\begin{eqnarray}\label{2.1}
{\mathbb E}[\sup_{0\leq t\leq 1}K_i^\gamma(t) ]\leq \sum_{(x,y)\in \Delta_i}{\mathbb E}[\sup_{0\leq t\leq 1}\|X_t(x)-X_t(y)\|_H^\gamma]\leq Cd2^{(i+1)d}2^{-id}\varphi(2^{-i})=C_1\varphi(2^{-i}).
\end{eqnarray}

For two points $x_1=(x_1^{(1)},x_1^{(2)},\ldots,x_1^{(d)})$ and $x_2=(x_2^{(1)},x_2^{(2)},\ldots,x_2^{(d)})$ in $[0,1]^d$, if $x_1^{(i)}\leq x_2^{(i)} (1\leq i\leq d)$, we write it as $x_1\leq x_2$. Let $(x,y)\in [0,1]^d\times [0,1]^d$. Then
there is a natural number $m_0\geq0$, and two increasing sequences $\{x_m\}$ and $\{y_m\} (x_m\in D_m,y_m\in D_m)$, such that for each $m<m_0$, $x_m\leq x$, $y_m\leq y$ and from each $m\geq m_0$, $x_m=x$, $y_m=y$.

If $|x-y|\leq 2^{-m}$, then either $x_m=y_m$ or $(x_m,y_m)\in \Delta_m$ and in any case
\begin{eqnarray}\label{2.2}
X_t(x)-X_t(y)=\sum_{i=m}^\infty(X_t(x_{i+1})-X_t(x_{i}))+X_t(x_m)-X_t(y_m)-
\sum_{i=m}^\infty(X_t(y_{i+1})-X_t(y_{i})).
\end{eqnarray}

Since points $x_{i+1}$ and $x_i$ can be connected by a piecewise line involving, for any $i\geq m$, at most $d$ steps between the nearest neighbors in $D_{i+1}$ such that each line segment's length is no more that $2^{-i-1}$, we conclude from (\ref{2.2}) that
\begin{eqnarray}\label{2.3}
\|X_t(x)-X_t(y)\|_H\leq K_m(t)+2d\sum_{i=m+1}^\infty K_i(t)\leq 2d\sum_{i=m}^\infty K_i(t).
\end{eqnarray}
Let $\alpha=1/\gamma-\vartheta$. Then
\begin{eqnarray}\label{2.4}
[X_t(\cdot)]_{\alpha,\varphi}&=&\sup_{x\neq y}\Big\{ \frac{\|X_t(x)-X_t(y)\|_H}{\varphi^\alpha(|x-y|)}, \ x,y\in [0,1]^d\Big\} \nonumber\\&\leq&
\sup_{x\neq y,m\in {\mathbb N}}\Big\{ \varphi^{-\alpha}(2^{-m-1})\sup_{2^{-m-1}\leq |x-y|< 2^{-m}} \|X_t(x)-X_t(y)\|_H, \ x,y\in [0,1]^d \Big\}
\nonumber\\&&+
\sup_{x\neq y}\Big\{ \varphi^{-\alpha}(1)\sup_{|x-y|\geq 1} \|X_t(x)-X_t(y)\|_H, \ x,y\in [0,1]^d \Big\}.
\end{eqnarray}

Observing that for any $(x,y)\in [0,1]^d\times [0,1]^d$ with $|x-y|\geq 1$, can be connected by a piecewise linear path  at most $d$ steps  such that every line segment's length is no more than $1$. Therefore
\begin{eqnarray}\label{2.5}
&&\sup_{x\neq y}\Big\{ \varphi^{-\alpha}(1)\sup_{|x-y|\geq 1} \|X_t(x)-X_t(y)\|_H, \ x,y\in [0,1]^d \Big\}
\nonumber\\&\leq & d\sup_{x\neq y}\Big\{ \varphi^{-\alpha}(1)\sup_{ |x-y|\leq1} \|X_t(x)-X_t(y)\|_H, \ x,y\in [0,1]^d \Big\}.
\end{eqnarray}
Combining (\ref{2.4}) and (\ref{2.5}), with the aid of (\ref{2.3}), we conclude that
\begin{eqnarray}\label{2.6}
[X_t(\cdot)]_{\alpha,\varphi}&\leq&
\sup_{x\neq y,m\in {\mathbb N}}\Big\{ \varphi^{-\alpha}(2^{-m-1})\sup_{ |x-y|\leq 2^{-m}} \|X_t(x)-X_t(y)\|_H, \ x,y\in [0,1]^d \Big\}
\nonumber\\&&+
d\sup_{x\neq y}\Big\{ \varphi^{-\alpha}(1)\sup_{|x-y|\leq 1} \|X_t(x)-X_t(y)\|_H, \ x,y\in [0,1]^d \Big\}
\nonumber\\&\leq &(d+1)
\sup_{x\neq y,m\in {\mathbb N}}\Big\{ \varphi^{-\alpha}(2^{-m-1})\sup_{ |x-y|\leq 2^{-m}} \|X_t(x)-X_t(y)\|_H, \ x,y\in [0,1]^d \Big\}
\nonumber\\&\leq&
2d(d+1)\sup_{x\neq y,m\in {\mathbb N}}\Big\{\varphi^{-\alpha}(2^{-m-1})\sum_{i=m}^\infty K_i(t) \Big\}\nonumber\\&\leq&
2d(d+1)\sup_{m\in\mN}\frac{ \varphi^{-\alpha}(2^{-m-1})}{\varphi^{-\alpha}(2^{-m})} \sup_{x\neq y,m\in {\mathbb N}}\Big\{\varphi^{-\alpha}(2^{-m})\sum_{i=m}^\infty K_i(t) \Big\}
\nonumber\\&\leq &
2d(d+1)\sup_{m\in\mN}\frac{ \varphi^{-\alpha}(2^{-m-1})}{\varphi^{-\alpha}(2^{-m})}\sum_{i=0}^\infty\varphi^{-\alpha}(2^{-i}) K_i(t).
\end{eqnarray}
Thanks to (\ref{1.5}), there is a constant $C>0$ such that
\begin{eqnarray*}
\sup_{m\in\mN}\frac{ \varphi^{-\alpha}(2^{-m-1})}{\varphi^{-\alpha}(2^{-m})} \leq C.
\end{eqnarray*}
Hence we achieve from  (\ref{2.6}) that
\begin{eqnarray}\label{2.7}
[X_t(\cdot)]_{\alpha,\varphi}\leq
2d(d+1)
C \sum_{i=0}^\infty\varphi^{-\alpha}(2^{-i}) K_i(t):=C(d)\sum_{i=0}^\infty\varphi^{-\alpha}(2^{-i}) K_i(t).
\end{eqnarray}

For $\gamma\geq 1$, by using Minkowski's inequality, we get from (\ref{2.1}) and (\ref{2.7}) that
\begin{eqnarray}\label{2.8}
[{\mathbb E} \sup_{0\leq t\leq 1}[X_t(\cdot)]_{\alpha,\varphi}^\gamma]^{\frac{1}{\gamma}} &\leq&
 C(d)\sum_{i=0}^\infty\varphi^{-\alpha}(2^{-i}) [{\mathbb E}\sup_{0\leq t\leq 1}K_i^\gamma(t)]^{\frac{1}{\gamma}}
\nonumber\\&\leq& C(d)C_1^{\frac{1}{\gamma}} \sum_{i=0}^\infty\varphi^{-\alpha}(2^{-i}) \varphi^{\frac{1}{\gamma}}(2^{-i})\nonumber\\&=&C(d)C_1^{\frac{1}{\gamma}} \sum_{i=0}^\infty\varphi^{\vartheta}(2^{-i}).
\end{eqnarray}

Since the series in (\ref{2.2}) are actually finite sums,  for $\gamma< 1$, we gain an analogue of (\ref{2.3})
\begin{eqnarray}\label{2.9}
\|X_t(x)-X_t(y)\|_H^\gamma\leq K_m^\gamma(t)+2d\sum_{i=m+1}^\infty K_i^\gamma(t)\leq 2d\sum_{i=m}^\infty K_i^\gamma(t).
\end{eqnarray}
The same calculations applied to $[X_t(\cdot)]_{\alpha,\varphi}$ in (\ref{2.4}) is adapted to $[X_t(\cdot)]_{\alpha,\varphi}^\gamma$ yields that: there is a constant $C(d)>0$ such that
\begin{eqnarray}\label{2.10}
[X_t(\cdot)]_{\alpha,\varphi}^\gamma\leq  C(d) \sum_{i=0}^\infty\varphi^{-\alpha\gamma}(2^{-i}) K_i^\gamma(t).
\end{eqnarray}
Combining (\ref{2.9}) and (\ref{2.10}), we arrive at
\begin{eqnarray}\label{2.11}
[{\mathbb E} \sup_{0\leq t\leq 1}[X_t(\cdot)]_{\alpha,\varphi}^\gamma] &\leq&
 C(d)\sum_{i=0}^\infty\varphi^{-\alpha\gamma}(2^{-i}) [{\mathbb E}\sup_{0\leq t\leq 1}K^\gamma_i(t)]
\nonumber\\&\leq& C(d)C_1 \sum_{i=0}^\infty\varphi^{\gamma\vartheta}(2^{-i}).
\end{eqnarray}

From (\ref{2.8}) and (\ref{2.11}), in view of (\ref{1.3}) and (\ref{1.4}), then
\begin{eqnarray}\label{2.12}
[{\mathbb E} \sup_{0\leq t\leq 1}[X_t(\cdot)]_{\alpha,\varphi}^\gamma]^{\frac{1}{\gamma}} <+\infty.
\end{eqnarray}
It follows from (\ref{2.12}) that for almost every $\omega$, $X_t(\cdot)$ is uniformly continuous on $[0,1]^d$ and it is uniformly in $t$, so it make sense to set
\begin{eqnarray*}
\tilde{X}_t(x,\omega)=\lim_{y\in [0,1]^d, y\rightarrow x}X_t(y,\omega).
\end{eqnarray*}
By Fatou's lemma and the hypothesis,
\begin{eqnarray*}
\mP\{\tilde{X}_t(x,\omega)=X_t(x,\omega), \ t\in [0,1]\}=1.
\end{eqnarray*}
Therefore, $\tilde{X}$ is the desired realization.

For general $\alpha\in (0,1/\gamma-\vartheta)$, we have
\begin{eqnarray}\label{2.13}
[X_t(\cdot)]_{\alpha,\varphi}&=&\sup_{x\neq y}\Big\{ \frac{\|X_t(x)-X_t(y)\|_H}{\varphi^\alpha(|x-y|)}, \ x,y\in [0,1]^d\Big\}\nonumber\\
&=&\sup_{x\neq y}\Big\{ \frac{\|X_t(x)-X_t(y)\|_H}{\varphi^{\frac1\gamma-\vartheta}(|x-y|)}\times
\frac{\varphi^{\frac1\gamma-\vartheta}(|x-y|)}{\varphi^{\alpha}(|x-y|)}, \ x,y\in [0,1]^d\Big\}.
\end{eqnarray}
Notice that $\varphi$ is continuous and nondecreasing, if $\varphi(\sqrt{d})\leq 1$, then
\begin{eqnarray*}
\sup_{x\neq y}\Big\{ \frac{\varphi^{\frac1\gamma-\vartheta}(|x-y|)}{\varphi^{\alpha}(|x-y|)}\ x,y\in [0,1]^d\Big\}\leq 1.
\end{eqnarray*}
If $\varphi(\sqrt{d})> 1$, then there is $0<\tilde{r}<\sqrt{d}$ such that $\varphi(\tilde{r})= 1$. Therefore,
\begin{eqnarray}\label{2.14}
\sup_{x\neq y}\Big\{ \frac{\varphi^{\frac1\gamma-\vartheta}(|x-y|)}{\varphi^{\alpha}(|x-y|)}\ x,y\in [0,1]^d\Big\}&\leq& \sup_{\tilde{r}\leq |x-y|\leq \sqrt{d}}\Big\{ \frac{\varphi^{\frac1\gamma-\vartheta}(|x-y|)}{\varphi^{\alpha}(|x-y|)}\Big\}+ 1\nonumber\\&\leq& \frac{\sup_{\tilde{r}\leq r\leq \sqrt{d}}\varphi^{\frac1\gamma-\vartheta}(r)}{\inf_{\tilde{r}\leq r\leq \sqrt{d}}\varphi^{\alpha}(r)}+ 1<+\infty.
\end{eqnarray}
In view of (\ref{2.13}) and (\ref{2.14}), (\ref{2.12}) is true, thus we complete the proof.
$\Box$

\section{Application to nonlocal stochastic heat equations}\label{sec}
\setcounter{equation}{0}
\begin{defi} \label{def3.1} (\cite{App}) Let $(\Omega,{\mathcal F},\{{\mathcal F}_t\}_{t\geq0},{\mathbb P})$ be a filtered complete probability space with the right continuous filtration ${\mathcal F}_t$.
Let $E$ be a ball $B_c(0)\backslash \{0\}$ in $\mathbb{R}^d$,  of radius $c$ without the center.
A time homogeneous Poisson random measure $N$ on $(E,\mathcal{B}(E))$ over the filtered probability space $(\Omega,{\mathcal F},\{{\mathcal F}_t\}_{t\geq0},{\mathbb P})$ with an intensity measure $\nu\times \lambda$, is a measurable function $N: (\Omega, {\mathcal F})\rightarrow(M_+(E\times{\mathbb R}_+), \mathcal{M}_+(E\times{\mathbb R}_+))$, such that

(i) for each $B\times I \in \mathcal{B}(E)\times  \mathcal{B}({\mathbb R}_+)$, if $\nu(B)<\infty$, $N(B\times I)$ is a Poisson random variable with parameter $\nu(B)\lambda(I)$;

(ii) $N$ is independently scattered, i.e. if the sets $E_j\times I_j\in \mathcal{B}(E)\times  \mathcal{B}({\mathbb R}_+), \ j=1,\ldots,n$ are pairwise disjoint, then the random variables $N(B_j\times I_j), \  j=1,\ldots,n$ are mutually independent;

(iii) for each $U\in \mathcal{B}(E)$, the $\overline{\mathbb{N}} \ (= {\mN} \cup\{+\infty\})$-valued process $\{N((0,t],U)\}_{t>0}$ is $\{{\mathcal F}_t\}_{t\geq0}$-adapted and its increments are independent of the past.
\end{defi}

\begin{rem}\label{rem3.1}
 In the above definition,  $M_+(E\times {\mathbb R}_+)$ is the family of all $\sigma$-finite positive measures on $E\times {\mathbb R}_+$, $\lambda$ is the Lebesgue measure on ${\mathbb R}_+$, $\nu$ is a L\'{e}vy measure which satisfies
$$
\int_E1\wedge |v|^2\nu(dv)<\infty.
$$
\end{rem}

\begin{defi} \label{def3.2} Let $N$ be a homogeneous Poisson random measure on $(E,\mathcal{B}(E))$ over the probability space $(\Omega,{\mathcal F},\{{\mathcal F}_t\}_{t\geq0},{\mathbb P})$. The ${\mathbb R}$-valued process $\{\tilde{N}((0,t],A)\}_{t>0}$ defined by
$$
\tilde{N}((0,t],A)=N((0,t],A)-\nu(A)t, \quad t>0, \ A\in \mathcal{B}(E),
$$
is called a compensator Poisson random measure. And now $\{\tilde{N}((0,t],A)\}_{t>0}$ is a martingale on $(\Omega,{\mathcal F},\{{\mathcal F}_t\}_{t\geq0},{\mathbb P})$.
\end{defi}

Let $\alpha\in (0,2]$ and let $g\in L^1(\Omega;L^1_{loc}([0,\infty);L^1(E,\nu;L^\infty({\mathbb R}^d))))$ such that $\{g(t,v,x,\cdot)\}_{t\geq 0}$ as a family of $L^1(\Omega,{\mathcal F},{\mathbb P})$-valued random variables is $\{{\mathcal F}_t\}_{t\geq0}$-adapted. Consider the following Cauchy problem:
\begin{eqnarray}\label{3.1}
\left\{\begin{array}{ll}
du(t,x)=\Delta^{\frac\alpha2} u(t,x)dt+\int_Eg(t,x,v)\tilde{N}(dt,dv), \ t>0, \ x\in {\mathbb R}^d, \\
u(0,x)=0, \  x\in{\mathbb R}^d,  \end{array}\right.
\end{eqnarray}
where $\Delta^{\frac\alpha2}:=-(-\Delta)^{\frac{\alpha}{2}}$ and $(-\Delta)^{\frac{\alpha}{2}}$ is the fractional Laplacian on ${\mathbb R}^d$, and for $\alpha\in (0,2)$ defined by 
\begin{eqnarray*}
(-\Delta_x)^{\frac{\alpha}{2}}\phi(x)=c(d,\alpha) \mbox{P.V.} \int_{\mR^d}\frac{\phi(x)-\phi(z+x)}
{|z|^{d+\alpha}}dz, \quad \forall \ \phi \in \cD(\mR^d), \ x\in \mR^d, 
\end{eqnarray*}
with  $c(d,\alpha)=\alpha2^{\alpha-1}\pi^{-d/2}\Gamma{(\frac{d+\alpha}{2})}/\Gamma{(\frac{2-\alpha}{2})}$.

We call $u$ a mild solution of (\ref{3.1}) if the measurable function given by
\begin{eqnarray}\label{3.2}
u(t,x)=\int_{(0,t]} \int_EK(t-r,\cdot)\ast g(r,\cdot,v)(x)\tilde{N}(dr,dv),
\end{eqnarray}
where $K(t,x)$ is the heat kernel of the equation
\begin{eqnarray*}
\partial_tu(t,x)=\Delta^{\frac\alpha2} u(t,x), \ t>0, \ x\in \mR^d,
\end{eqnarray*}
satisfies:

(1) $\{u(t)\}_{t\geq0}$ is $\{{\mathcal F}_t\}_{t\geq0}$-adapted;

(2) $\{u(t,x,\cdot)\}_{t\geq 0}$ as a family of $L^1(\Omega,{\mathcal F},{\mathbb P})$-valued random variables is right continuous and has left limits in the variable $t>0$, namely,
\begin{eqnarray*}
u(t-,x,\cdot)=L^1(\Omega)-\lim_{s\uparrow t}u(s,x,\cdot), \ t>0.
\end{eqnarray*}

From the obvious representation (\ref{3.2}), by analogue calculations in \cite{AWZ,ZBL}, $u$ satisfies the properties claimed above. Now we will show that if $g$ is more regular, the stochastic field $u$ will admit a continuous realization which is continuous in $x$.

\begin{theorem} \label{the3.1} Let $g\in L^1(\Omega;L^1_{loc}([0,\infty);L^1(E,\nu;L^\infty({\mathbb R}^d))))$ such that $\{g(t,v,x,\cdot)\}_{t\geq 0}$ as a family of $L^1(\Omega,{\mathcal F},{\mathbb P})$-valued random variables is $\{{\mathcal F}_t\}_{t\geq0}$-adapted.  Suppose that $p\geq 1$ and $0<\vartheta<1/p$. Let $\varphi$ be a nonnegative and nondecreasing continuous function on
${\mathbb R}_+=(0,+\infty)$ such that $\lim_{r\rightarrow 0+}\varphi(r)=0$, and (\ref{1.3}) and (\ref{1.5}) hold for $\varphi$. For every $x,y\in \mR^d$, we assume there is
a nonnegative measurable function $h\in L^p(\Omega;L^1_{loc}([0,\infty);L^1(E,\nu))\cap  L^p_{loc}([0,\infty);L^p(E,\nu)))$ such that
\begin{eqnarray*}
 |g(t,v,x)-g(t,v,y)|\leq h(t,v)|x-y|^\frac{d}{p}\varphi^{\frac{1}{p}}(|x-y|).
\end{eqnarray*}
Let $u$ be defined by (\ref{3.2}). Then  $u$ admits a realization $\tilde{u}$ which is continuous in $x$, and for every $\beta\in (0,1/p-\vartheta]$, every $T>0$ and every $c_1>0$,
\begin{eqnarray}\label{3.3}
{\mathbb E}\sup_{0\leq t\leq T}\sup_{x\neq y}\Big\{ \frac{|\tilde{u}(t,x)-\tilde{u}(t,y)|^p}{\varphi^{\beta p}(|x-y|)}, \ x,y\in B_{c_1}\Big\}<+\infty.
\end{eqnarray}
\end{theorem}

To prove Theorem \ref{the3.1}, we need a lemma.
\begin{lemma}\label{lem3.1}  Let $E=B_c(0)\setminus\{0\}$ ($0<c\in{\mathbb R}$) and $p\geq1$. Suppose $\psi\in L^p(\Omega;L^p_{loc}([0,\infty);L^p(E,\nu)))$
is an $\{{\mathcal F}_t\}_{t\geq 0}$-adapted stochastic process and
\begin{eqnarray*}
I_t=\int_{(0,t]}\int_E\psi(r,v)\tilde{N}(dr,dv).
\end{eqnarray*}
(i) (Kunita's first inequality \cite[Theorem 4.4.23]{App}) If $p\geq 2$ and $\psi\in L^p(\Omega;L^2_{loc}([0,\infty);L^2(E,\nu)))$ further. Then for every $t>0$, there exists a positive constant $C(p)>0$, such that
\begin{eqnarray*}
\mE[\sup_{0\leq s\leq t}|I_s|^{p}]\leq C(p)\Big\{\mE
\Big[\int^t_0\int_E|\psi(r,v)|^2\nu(dv)dr\Big]^{\frac{p}{2}}+
\mE\int^t_0\int_E|\psi(r,v)|^p\nu(dv)dr \Big\}.
\end{eqnarray*}

(ii) (\cite[Proposition 2.2]{ZBL}) If $1\leq p<2$, then for every $t>0$, there exists a positive constant $C(p)>0$, such that
\begin{eqnarray*}
\mE[\sup_{0\leq s\leq t}|I_s|^{p}]\leq C(p)
\mE\int^t_0\int_E|\psi(r,v)|^p\nu(dv)dr.
\end{eqnarray*}
\end{lemma}

\vskip2mm\noindent
\textbf{Proof of Theorem \ref{the3.1}.} Let $p$ be given in Theorem \ref{the3.1}. For every $x,y\in B_{c_1}$ ($c_1>0$ is any given real number), and every $T>0$, from (\ref{3.2}) we have
\begin{eqnarray}\label{3.4}
&&\mE \sup_{0\leq t\leq T}|u(t,x)-u(t,y)|^p \nonumber \\ &=&\mE \sup_{0\leq t\leq T}\Big|\int_{(0,t]} \int_E[K(t-r,\cdot)\ast g(r,v,\cdot)(x)-K(t-r,\cdot)\ast g(r,v,\cdot)(y)]\tilde{N}(dr,dv)\Big|^p.
\end{eqnarray}
Since $\tilde{N}((0,t],\cdot)=N((0,t],A)-\nu(\cdot)t$, from (\ref{3.4}), it yields that
\begin{eqnarray}\label{3.5}
&&\mE \sup_{0\leq t\leq T}|u(t,x)-u(t,y)|^p \nonumber \\ &=&\mE \sup_{0\leq t\leq T}\Big|\int_{(0,t]} \int_E[K(t-r,\cdot)\ast g(r,v,\cdot)(x)-K(t-r,\cdot)\ast g(r,v,\cdot)(y)]N(dr,dv)   \nonumber \\ && \quad -\int_{(0,t]} \int_E[K(t-r,\cdot)\ast g(r,v,\cdot)(x)-K(t-r,\cdot)\ast g(r,v,\cdot)(y)]\nu(dv)dr\Big|^p
   \nonumber \\
&\leq& 2^{p-1}\mE \sup_{0\leq t\leq T}\Big|\int_{(0,t]} \int_E|K(t-r,\cdot)\ast g(r,v,\cdot)(x)-K(t-r,\cdot)\ast g(r,v,\cdot)(y)|N(dr,dv)\Big|^p   \nonumber \\ && + 2^{p-1}\mE \sup_{0\leq t\leq T}\Big|\int_{(0,t]} \int_E|K(t-r,\cdot)\ast g(r,v,\cdot)(x)-K(t-r,\cdot)\ast g(r,v,\cdot)(y)|\nu(dv)dr\Big|^p.
\end{eqnarray}

We estimate the convolution in (\ref{3.5}) by
\begin{eqnarray}\label{3.6}
&&|K(t-r,\cdot)\ast g(r,v,\cdot)(x)-K(t-r,\cdot)\ast g(r,v,\cdot)(y)|
\nonumber \\ &=&  \Big|\int_{\mR^d}K(t-r,z)[g(r,v,x-z)-g(r,v,y-z)]dz\Big|
\nonumber \\ &\leq& \int_{\mR^d}K(t-r,z)|g(r,v,x-z)-g(r,v,y-z)|dz
\nonumber \\ &\leq& h(r,v)|x-y|^{\frac{d}{p}}\varphi^{\frac{1}{p}}(|x-y|).
\end{eqnarray}

Combining (\ref{3.5}) and (\ref{3.6}), we get
\begin{eqnarray}\label{3.7}
&&\mE \sup_{0\leq t\leq T}|u(t,x)-u(t,y)|^p \nonumber \\
&\leq& C\mE \sup_{0\leq t\leq T}\Big|\int_{(0,t]} \int_Eh(r,v)N(dr,dv)\Big|^p|x-y|^{d}\varphi(|x-y|) \nonumber \\ && + C\mE\Big|\int_0^T \int_Eh(r,v)\nu(dv)dr\Big|^p|x-y|^{d}\varphi(|x-y|)
\nonumber \\
&\leq& C\mE \sup_{0\leq t\leq T}\Big|\int_{(0,t]} \int_Eh(r,v)\tilde{N}(dr,dv)\Big|^p|x-y|^{d}\varphi(|x-y|) \nonumber \\ && + C \mE\Big|\int_{(0,T]} \int_Eh(r,v)\nu(dv)dr\Big|^p|x-y|^{d}\varphi(|x-y|) .
\end{eqnarray}
If $p\geq2$, by using Lemma \ref{lem3.1} (i), we conclude from (\ref{3.7}) that
\begin{eqnarray}\label{3.8}
&&\mE \sup_{0\leq t\leq T}|u(t,x)-u(t,y)|^p \nonumber \\
&\leq& C\mE\Big|\int_0^T \int_E|h(r,v)|^2\nu(dv)dr\Big|^{\frac{p}{2}}|x-y|^{d}\varphi(|x-y|)
\nonumber \\ && +
C\mE\int_0^T\int_E|h(r,v)|^p\nu(dv)dr|x-y|^{d}\varphi(|x-y|)
\nonumber \\ && + C \mE\Big|\int_0^T \int_Eh(r,v)\nu(dv)dr\Big|^p|x-y|^{d}\varphi(|x-y|).
\end{eqnarray}
Observing that $h\in L^p(\Omega; L^1_{loc}([0,\infty);L^1(E,\nu))\cap  L^p_{loc}([0,\infty);L^p(E,\nu)))$, therefore
\begin{eqnarray}\label{3.9}
\mE \sup_{0\leq t\leq T}|u(t,x)-u(t,y)|^p
\leq C|x-y|^{d}\varphi(|x-y|).
\end{eqnarray}

If $1\leq p<2$, with the aid Lemma \ref{lem3.1} (ii), we conclude from (\ref{3.7}) that
\begin{eqnarray*}
&&\mE \sup_{0\leq t\leq T}|u(t,x)-u(t,y)|^p \nonumber \\
&\leq&
C\mE\int_0^T\int_E|h(r,v)|^p\nu(dv)dr|x-y|^{d}\varphi(|x-y|)
\nonumber \\ && + C \mE\Big|\int_0^T \int_Eh(r,v)\nu(dv)dr\Big|^p|x-y|^{d}\varphi(|x-y|),
\end{eqnarray*}
which also implies the estimate (\ref{3.9}). In view of Remark \ref{rem1.1} (i), we finish the proof. $\Box$

\begin{rem}\label{rem3.2} The main ingredient in proving Theorem \ref{the3.1} is to estimate the stochastic evolution
\begin{eqnarray}\label{3.10}
\mE \sup_{0\leq t\leq T}\Big\|\int_{(0,t]} \int_E[f(r,v,x) \tilde{N}(dr,dv)\Big\|_{L^\infty_x(\mR^d)}^p
\end{eqnarray}
for $f\in L^p(\Omega;L^1_{loc}([0,\infty);L^1(E,\nu;L^\infty(\mR^d)))\cap  L^p_{loc}([0,\infty);L^p(E,\nu;L^\infty(\mR^d))))$. When the integrand and $L^\infty(\mR^d)$ are replaced by $U(t,r)h(r-)$ ($U$ is a evolution operator) and a Hilbert space, respectively, the estimate for (\ref{3.10}) was established by Kotelenez \cite{Kot84}. Since then, Kotelenez's result was strengthened by Hamedani and Zangeneh \cite{HZ}, Ichikawa \cite{Ich}. We refer to  \cite{BH,DMN,Hau,HS08,MR,WDL} for some other extensions. Observing that all extensions are concentrated on stochastic evolution taking values in martingale type ($1<p<\infty$) Banach spaces, and as noticed in  \cite[Remark 2.11]{Hau}, $L^\infty(\mR^d)$ is not a Banach space of martingale type $p$ for any $p>1$, so the following estimates:
\begin{eqnarray*}
&&\mE \sup_{0\leq t\leq T}\Big\|\int_{(0,t]} \int_Ef(r,v,x) \tilde{N}(dr,dv)\Big\|_{L^\infty_x(\mR^d)}^p
\\ &\leq& C(p)\Big\{\mE
\Big[\int^T_0\int_E\|f(r,v,\cdot)\|^2_{L^\infty(\mR^d)}\nu(dv)dr\Big]^{\frac{p}{2}} \\ && +
\mE\int^T_0\int_E\|f(r,v,\cdot)\|^p_{L^\infty(\mR^d)}\nu(dv)dr \Big\}, \ \ p\geq 2,
\end{eqnarray*}
and
\begin{eqnarray*}
&&\mE \sup_{0\leq t\leq T}\Big\|\int_{(0,t]} \int_E[f(r,v,x) \tilde{N}(dr,dv)\Big\|_{L^\infty_x(\mR^d)}^p
\\ &\leq& C(p)
\mE\int^T_0\int_E\|f(r,v,\cdot)\|^p_{L^\infty(\mR^d)}\nu(dv)dr, \ \ 1< p< 2,
\end{eqnarray*}
for Poisson random measure in general will be not true. Fortunately, if we assume further that $f\in L^p(\Omega;L^1_{loc}([0,\infty);L^1(E,\nu;L^\infty(\mR^d))))$, then we gain
\begin{eqnarray}\label{3.11}
&&\mE \sup_{0\leq t\leq T}\Big\|\int_{(0,t]} \int_E[f(r,v,x) \tilde{N}(dr,dv)\Big\|_{L^\infty_x(\mR^d)}^p
\nonumber\\ &\leq& C(p)\Big\{\mE
\Big[\int^T_0\int_E\|f(r,v,\cdot)\|_{L^\infty(\mR^d)}\nu(dv)dr\Big]^{p}\nonumber\\&& +
\mE\int^T_0\int_E\|f(r,v,\cdot)\|^p_{L^\infty(\mR^d)}\nu(dv)dr \Big\}.
\end{eqnarray}
The above estimate will play  a central role in stochastic partial equations when we prove the boundedness of solutions.
\end{rem}

From (\ref{3.11}), for the Cauchy problem (\ref{3.1}), we have the following corollary.

\begin{corollary} \label{cor3.1} Let $g\in L^1(\Omega;L^1_{loc}([0,\infty);L^1(E,\nu;L^\infty({\mathbb R}^d))))$ such that $\{g(t,v,x,\cdot)\}_{t\geq 0}$ as a family of $L^1(\Omega,{\mathcal F},{\mathbb P})$-valued random variables is $\{{\mathcal F}_t\}_{t\geq0}$-adapted, and let $u$ be given by (\ref{3.2}). Then $u\in L^1(\Omega; L^\infty_{loc}([0,\infty);L^1(E,\nu;L^\infty(\mR^d))))$.
\end{corollary}
\vskip2mm\noindent
\textbf{Proof.} Let $T>0$, from (\ref{3.2}) we have
\begin{eqnarray}\label{3.12}
&&\mE \sup_{0\leq t\leq T}ess \sup_{x\in\mR^d}|u(t,x)|\nonumber \\ &=&\mE \sup_{0\leq t\leq T}ess \sup_{x\in\mR^d}\Big|\int_{(0,t]} \int_EK(t-r,\cdot)\ast g(r,v,\cdot)(x)\tilde{N}(dr,dv)\Big|.
\nonumber \\
&\leq& \mE \sup_{0\leq t\leq T}ess \sup_{x\in\mR^d}\Big|\int_{(0,t]} \int_E|K(t-r,\cdot)\ast g(r,v,\cdot)(x)|N(dr,dv)\Big|  \nonumber \\ && +\mE \sup_{0\leq t\leq T}ess \sup_{x\in\mR^d}\Big|\int_{(0,t]} \int_E|K(t-r,\cdot)\ast g(r,v,\cdot)(x)|\nu(dv)dr\Big|
\nonumber \\
&\leq& \mE \sup_{0\leq t\leq T}\Big|\int_{(0,t]} \int_E\| g(r,v)\|_{L^\infty_x(\mR^d)}N(dr,dv)\Big|  \nonumber \\ && +\mE \sup_{0\leq t\leq T}\Big|\int_{(0,t]} \int_E\| g(r,v)\|_{L^\infty_x(\mR^d)}\nu(dv)dr\Big|
\nonumber \\
&\leq& \mE \sup_{0\leq t\leq T}\Big|\int_{(0,t]} \int_E\| g(r,v)\|_{L^\infty_x(\mR^d)}\tilde{N}(dr,dv)\Big|  \nonumber \\ && +2\mE \sup_{0\leq t\leq T}\Big|\int_{(0,t]} \int_E\| g(r,v)\|_{L^\infty_x(\mR^d)}\nu(dv)dr\Big| \nonumber \\
&\leq& C\mE\int_0^T\int_E\| g(r,v)\|_{L^\infty_x(\mR^d)}\nu(dv)dr<+\infty.
\end{eqnarray}

\begin{rem}\label{rem3.3} (i) Let $q>1$ and  let $u$ be given by (\ref{3.2}). If $g\in L^1(\Omega;L^1_{loc}([0,\infty);L^1(E,\nu;L^q({\mathbb R}^d))))$, we also prove that $u\in L^1(\Omega; L^\infty_{loc}([0,\infty);L^1(E,\nu;L^q(\mR^d))))$.

(ii) For regularities of solutions to equation (\ref{3.2}), one also consults to
\cite{KK2016,KK2012,LGWW,MPR}.
\end{rem}

\vskip2mm\noindent {\bf Acknowledgment.} This work was supported in part
by NSFC of China grants 11501577, 11771123.

 \end{document}